\documentclass{amsart}
\usepackage{amssymb,latexsym,epsfig,color}

\newcommand{\Z}{{\mathbb Z}}
\newcommand{\N}{{\mathbb N}}
\newcommand{\R}{{\mathbb R}}

\newcommand{\PR}{{\mathbb P}}

\newcommand{\A}{{\mathcal A}}

\newcommand{\B}{{\mathcal B}}
\newcommand{\LL}{{\mathcal L}}
\newcommand{\CC}{{\mathcal C}}

\newtheorem{theorem}{Theorem}[section]

\newtheorem{corollary}[theorem]{Corollary}
\newtheorem{conjecture}[theorem]{Conjecture}
\newtheorem{lemma}[theorem]{Lemma}

\newtheorem{question}[theorem]{Question}

\theoremstyle{definition}
\newtheorem{definition}[theorem]{Definition}
\newtheorem{remark}[theorem]{Remark}
\newtheorem{example}[theorem]{Example}
\newtheorem{algorithm}[theorem]{Algorithm}

\title [$\Delta$-Normal Configurations]
{Toric Initial Ideals of $\Delta$-Normal Configurations:
\newline Cohen-Macaulayness and Degree Bounds}

\author{Edwin O'Shea}
\author{Rekha R. Thomas}
\thanks{Research partially supported by NSF grant DMS 0100141.
The first author is partially supported by a Fulbright grant.}
\address{Department of Mathematics, University of
  Washington, Seattle, WA 98195-4350}
\email{[oshea,thomas]@math.washington.edu}
\date{\today}

\begin{document}

\begin{abstract} 
  A normal (respectively, graded normal) vector configuration $\A$
  defines the toric ideal $I_{\mathcal A}$ of a normal (respectively,
  projectively normal) toric variety. These ideals are Cohen-Macaulay,
  and when $\mathcal A$ is normal and graded, $I_{\mathcal A}$ is
  generated in degree at most the dimension of $I_{\A}$. Based on
  this, Sturmfels asked if these properties extend to initial ideals
  --- when $\mathcal A$ is normal, is there an initial ideal of
  $I_{\mathcal A}$ that is Cohen-Macaulay, and when $\mathcal A$ is
  normal and graded, does $I_{\mathcal A}$ have a Gr\"obner basis
  generated in degree at most $dim(I_{\A})$~?  In this paper, we
  answer both questions positively for $\Delta$-normal configurations.
  These are normal configurations that admit a regular triangulation
  $\Delta$ with the property that the subconfiguration in each cell of
  the triangulation is again normal.  Such configurations properly
  contain among them all vector configurations that admit a regular
  unimodular triangulation. We construct non-trivial families
  of both $\Delta$-normal and non-$\Delta$-normal configurations.
\end{abstract}

\maketitle
\section{Introduction}

A finite vector configuration $\A = \left\{ {\bf a}_i \, : \, i = 1,
  \ldots, n \right\} \subset \Z^d$ defines the {\em toric ideal} $
I_\A := \langle {\bf x^u} - {\bf x^v} \, : \, {\bf u}, {\bf v} \in
\N^n , \,\, \sum_{i = 1}^{n} {\bf a}_i u_i = \sum_{i=1}^{n} {\bf a}_i
v_i \rangle $ in the polynomial ring $R := {\mathbb K} [x_1,
\ldots,x_n] = {\mathbb K} [{\bf x}]$ where $\mathbb K$ is a field. Let
$cone(\A), \,\, \Z \A$ and $\N \A$ denote the cone, lattice and
semigroup spanned by the $\R_{\geq 0}, \Z$ and $\N$-linear
combinations of $\A$ where $\N$ is the set of non-negative integers.
Let $dim(\A)$ be the Krull dimension of $R/I_{\A}$ which equals the
rank of $\Z \A$.  Assume $dim(\A) = d$.  The configuration $\A$ is
{\em normal} if $\N \A = cone(\A) \cap \Z \A$ and {\em graded} if $\A$
spans an affine hyperplane in $\R^d$.  A finite set $\B \subset \Z^d$
such that $\N \B = cone(\A) \cap \Z \A$ is called a {\em Hilbert
  basis} of the semigroup $cone(\A)\cap \Z \A $.  If $\A$ is normal,
the zero set of $I_{\A}$ is a normal toric variety in ${\mathbb K}^n$
of dimension $d$, and when $\A$ is also graded, it is a projectively
normal toric variety in $\PR^{n-1}_{\mathbb K}$ of dimension $d-1$.
See \cite{stu} for details on toric ideals. A survey of recent results
and open questions on normal configurations can be found in
\cite{bgt2}.

It is well known that initial ideals of a polynomial ideal inherit
important invariants of the ideal such as dimension, degree and
Hilbert function. Thus it is natural to ask if certain initial ideals
inherit further properties of the ideal such as Cohen-Macaulayness,
Betti numbers or reducedness (of the corresponding scheme). A result
of Hochster~\cite{hochster} shows that when $\A$ is normal, $I_\A$ is
Cohen-Macaulay. If $\A$ is also graded, then $I_\A$ is generated by
homogeneous binomials of degree at most $d$ \cite[Thm.~13.14]{stu}.
Motivated by this, Sturmfels asked and conjectured the following.

\begin{question} \label{cmconjecture}
  If $\A$ is normal (more generally, if $I_{\A}$ is
  Cohen-Macaulay), is there a monomial initial ideal
  $in_{\omega}(I_\A)$ of the toric ideal $I_{\A}$ that is
  Cohen-Macaulay~?
\end{question}

\begin{conjecture} \label{bigconjecture} ~\cite[Conjecture 2.8]{etv}
  If $\A$ is a graded, normal configuration then $I_\A$ has a
  Gr\"obner basis whose elements have degree at most $d = dim(\A)$.
\end{conjecture}

In this paper, we show that Question~\ref{cmconjecture} has a positive
answer and Conjecture~\ref{bigconjecture} is true when $\A$ is
$\Delta$-{\em normal}. These configurations were defined by Ho{\c
  s}ten and Thomas in \cite{gomory}. We recall the definition. Let
$\Delta$ be a pure $(d-1)$-dimensional simplicial complex with vertex
set $[n] := \{1, \ldots, n\}$. We denote the set of {\em facets}
($d$-element faces) of $\Delta$ by $\textup{max} \, \Delta$.  For a
set $\tau \subseteq [n]$, let $\A_{\tau} := \{ {\bf a}_i \in \A \, :
\, i \in \tau \}$. We say that $\Delta$ is a {\em triangulation} of
$\A$ if $cone(\A) = \bigcup_{\sigma \in \textup{max} \, \Delta}
cone(\A_{\sigma})$ and $cone(\A_{\sigma_i}) \cap cone(\A_{\sigma_j}) =
cone(\A_{{\sigma_i} \cap {\sigma_j}})$ for all $\sigma_i, \sigma_j \in
\textup{max} \, \Delta$. The {\em Stanley-Reisner} ideal of $\Delta$
is the squarefree monomial ideal $\langle \prod_{i \in \tau} x_i \, :
\, \tau \not \in \Delta \rangle \subseteq R$.  A cornerstone in the
combinatorics of toric initial ideals is the result by Sturmfels that
the radical of a monomial initial ideal $in_{\omega}(I_{\A})$ of the
toric ideal $I_{\A}$ is the Stanley-Reisner ideal of a triangulation
$\Delta_{\omega}$ of $\A$ \cite[Thm.~8.3]{stu}. The ideal
$in_{\omega}(I_{\A})$ is said to be {\em supported on}
$\Delta_{\omega}$. Triangulations supporting initial ideals of
$I_{\A}$ are the {\em regular triangulations} of $\A$.  A
triangulation $T$ of $\A$ is {\em unimodular} if for each $\sigma \in
\textup{max} \, T$, $\Z \A_{\sigma} = \Z \A$.

\begin{definition} \label{def_dnormal}
A configuration $\A$ is {\bf $\Delta$-normal} if it has a regular
triangulation $\Delta$ such that for each $\sigma \in \textup{max} \,
\Delta$, $\A \cap cone(\A_{\sigma})$ is a Hilbert basis of
$cone(\A_{\sigma}) \cap \Z \A$. 
\end{definition}

Note that $\Z \A$ is used in the semigroups of
Definition~\ref{def_dnormal}. All $\Delta$-normal configurations are 
normal. In Sections~\ref{cmsection} and \ref{dnormalboundsection}, we
prove our main results.

\vspace{.2cm} 

\noindent {\bf Theorem~\ref{cm}.} Let $\A$ be a
$\Delta$-normal configuration. Then there exists a term order $\succ$
such that $\Delta = \Delta_\succ$ and $in_{\succ}(I_\A)$ is
Cohen-Macaulay.

\vspace{.2cm} 

\noindent {\bf Theorem~\ref{dnormalbound}.} Let $\A$ be a
graded $\Delta$-normal configuration. Then there exists a term order
$\succ$ such that $\Delta = \Delta_\succ$ and the Gr\"obner basis of 
$I_{\A}$ with respect to $\succ$ consists of binomials of degree at
most $d = dim(\A)$. 

\vspace{.2cm} 

It was shown in~\cite{gomory} that if $\A$ is $\Delta$-normal then
$I_{\A}$ has a monomial initial ideal that is free of embedded primes.
In Section~\ref{setup} we recall the main features of this initial
ideal. Theorems~\ref{cm} and \ref{dnormalbound} are proved by showing
that this same initial ideal is Cohen-Macaulay and, when $\A$ is
graded, generated in degree at most $d$.  Our proofs are combinatorial
and rely heavily on the structure of this initial ideal.

The set of $\Delta$-normal configurations is a proper subset of the
set of normal configurations.  They occur naturally in two ways.  If
$\A$ has a regular unimodular triangulation $\Delta$, then $\A$ is
$\Delta$-normal. If $cone(\A)$ is simplicial and we assume that its
extreme rays are generated by ${\bf a}_1, \ldots, {\bf a}_d$, then
$\A$ is $\Delta$-normal with respect to the coarsest regular
triangulation $\Delta$ consisting of the unique facet $\sigma = \{1,
\ldots, d\}$.  These were the only examples known so far.  Specific
instances of normal configurations that are not $\Delta$-normal for
any $\Delta$ are also known \cite{gomory}. In Section~\ref{families}
we construct non-trivial families of both $\Delta$-normal and
non-$\Delta$-normal configurations.

\vspace{.2cm} 

\noindent {\bf Theorem~\ref{babysfirstdeltanormalfamily}.}
There are families of $\Delta$-normal configurations $\{\A^d \subset
\Z^d, \, d \geq 5 \}$ where $cone(\A^d)$ is non-simplicial and
$\A^d$ has no regular unimodular triangulations.

\vspace{.2cm} 

\noindent {\bf Theorem~\ref{babysfirstnormalfamily}.} There is a 
family of normal, graded configurations $\{\A^d \subset \Z^d, \, d
\geq 11 \}$, that are not $\Delta$-normal for any regular
triangulation $\Delta$.

\section{Background: An initial ideal without embedded primes }
\label{setup}

We now recall from \cite{gomory} the initial ideal of $I_{\A}$ without
embedded primes when $\A$ is $\Delta$-normal. The construction uses
the {\em standard pair decomposition} of a monomial ideal $M$
\cite{stv} which carries detailed information about $Ass(M)$, the set
of associated primes of $M$.  Recall that every monomial prime ideal
of $R$ is of the form $P_{\tau} := \langle x_j \, : \, j \not \in \tau
\rangle$ for some $\tau \subseteq [n]$. The monomials of $R$ that do
not lie in $M$ are the {\em standard monomials} of $M$. The {\em
  support} of a monomial ${\bf x^v}$ is defined to be the support of
its exponent vector ${\bf v}$ --- i.e., $supp({\bf x^v}) = supp({\bf
  v}) := \{ i \in [n] \, : \, v_i \neq 0\}$.

\begin{definition}\cite{stv}
  Let $M \subseteq R$ be a monomial ideal. For a standard monomial
  ${\bf x^u}$ of $M$ and a set $\tau \subseteq [n]$, $({\bf x^u},
  \tau)$ is a {\bf pair} of $M$ if all monomials in ${\bf x^u} \cdot
  {\mathbb K}[x_j: j \in \tau]$ are standard monomials of $M$. We call
  $({\bf x^u}, \tau)$ a {\bf standard pair} of $M$ if:
\begin{enumerate}
\item$({\bf x^u}, \tau)$ is a pair of $M$,
\item $\tau \cap supp({\bf x^u}) = \emptyset$, and 
\item the set of monomials in ${\bf x^u} \cdot {\mathbb K}[x_j: j \in
  \tau]$ is not properly contained in the set of monomials in ${\bf
    x^v} \cdot {\mathbb K}[x_j: j \in \tau']$ for any $({\bf x^v},
  \tau')$ satisfying (1) and (2).
\end{enumerate}
\end{definition}

The set of standard pairs of $M$ is unique and is called the standard 
pair decomposition of $M$ since this set provides a decomposition of 
the standard monomials of $M$. If ${\bf x^v}$ is a
standard monomial of $M$ then there is a standard pair $({\bf x^u},
\tau)$ of $M$ such that ${\bf x^u}$ divides ${\bf x^v}$ and $supp({\bf
  x^{v-u}}) \subseteq \tau$.  In this case we say that ${\bf x^v}$ is
{\em covered} by $({\bf x^u}, \tau)$. We also use $({\bf x^u}, \tau)$
to denote the set of all standard monomials covered by it.

\begin{theorem} \label{ap56}
(A) \cite{stv} Let $M$ be a monomial ideal in $R$. Then,
\begin{enumerate} 
\item $P_{\tau} \in Ass(M)$ if and only if $M$ has a standard pair
  of the form $(\ast, \tau)$.
\item $P_{\tau}$ is a minimal prime of $M$ if and only if $(1, \tau)$
  is a standard pair of $M$.
\end{enumerate}
(B) \cite{stu} Let $M = in_{\omega}(I_{\A})$ be a monomial initial ideal
of the toric ideal $I_{\A}$. Then,
\begin{enumerate}
\item if $P_{\tau} \in Ass(M)$ then $\tau$ is a face of the regular 
triangulation $\Delta_{\omega}$ of $\A$, 
\item $P_{\sigma}$ is a minimal prime of $M$ if and only if $\sigma \in 
\textup{max} \, \Delta_{\omega}$, and 
\item for $\sigma \in \textup{max} \, \Delta_{\omega}$, the number of
  standard pairs of $M$ of the form $(\ast,\sigma)$ is $vol(\sigma)$,
  the normalized volume of $\sigma$ in $\Delta_{\omega}$.
\end{enumerate}
\end{theorem}

We call ${\bf x^u}$ and $\tau$ the {\em root} and {\em face} of the
standard pair $({\bf x^u}, \tau)$. Let ${\bf A}$ (respectively, ${\bf
  A}_{\sigma}$) be the matrix whose set of columns is $\A$
(respectively, $\A_{\sigma}$).  The normalized volume of $\sigma \in
\textup{max} \, \Delta_{w}$ is the absolute value of the determinant
of ${\bf A}_{\sigma}$ divided by the g.c.d. of the non-zero maximal
minors of ${\bf A}$.

\begin{theorem} \label{dnormal} \cite[Thm.~4.7]{gomory} Let $\A$ be
  a $\Delta$-normal configuration. Then there exists a term order
  $\succ$ such that $\Delta = \Delta_\succ$ and $in_\succ(I_\A)$ is
  free of embedded primes.
\end{theorem}

The term order $\succ$ needed in Theorem~\ref{dnormal} is described in
\cite{gomory} and is not directly used in this paper. The ideal
$in_\succ(I_\A)$ is shown to be free of embedded primes via an
explicit description of its standard pairs.  This structure is crucial
for this paper and hence we recall it now.  Assume without loss of
generality that $\Z \A = \Z^d$. For $\sigma \in \textup{max} \,
\Delta_{\succ}$, let $ FP_{\sigma} := \{\sum_{i \in \sigma}{\lambda_i
  {\bf a}_i } : 0 \leq \lambda_i < 1 \, \} $ be the half open {\em
  fundamental parallelopiped} of $cone(\A_{\sigma})$.  Then $FP_\sigma
\cap \Z ^d$ has $vol(\sigma)$ elements including the origin. For ${\bf
  \gamma} \in FP_\sigma \cap \Z ^d$, let ${\bf x^{u_{\gamma}}}$ be the
cheapest monomial with respect to $\succ$ among all ${\bf x^u} \in R$
with ${\bf A}{\bf u} = {\bf \gamma}$.  It was shown in \cite{gomory}
that $supp({\bf x^{u_{\gamma}}}) \subseteq \sigma_{in} := \{ i \, : \,
{\bf a}_i \in cone(\A_{\sigma}), i \notin \sigma \}$. The standard
pairs of the initial ideal $in_{\succ}(I_{\A})$ in
Theorem~\ref{dnormal} are precisely the pairs $({\bf x^{u_{\gamma}}},
\sigma)$ as ${\bf \gamma}$ varies in $FP_\sigma \cap \Z ^d$ for each
$\sigma$ in $\textup{max} \, \Delta_{\succ}$.  By Theorem~\ref{ap56},
$in_{\succ}(I_{\A})$ is thus free of embedded primes.

For the remainder of this paper we will denote the toric initial ideal
$in_{\succ}(I_\A)$ of Theorem~\ref{dnormal} by $J$ and its set of standard  
pairs by $\mathcal S(J)$. 

\begin{example} \label{example1}
Let $\A$ be the vector configuration consisting of the $13$ columns of
$${\bf A} = \left( \begin{array}{ccccccccccccc}
1 & 1 & 1 & 1 & 1 & 1 & 1 & 1 & 1 & 1 & 1 & 1 & 1\\
0 & 1 & 2 & 3 & 0 & 1 & 2 & 3 & 0 & 1 & 2 & 3 & 0\\
0 & 0 & 0 & 0 & 1 & 1 & 1 & 1 & 2 & 2 & 2 & 2 & 3 
\end{array} \right ).$$
Then $\A$ is a graded {\em supernormal} configuration \cite{HMS}  
which means that it is $\Delta$-normal with respect to every 
regular triangulation. Consider the regular triangulation  
$$
\Delta = \{ \{1,4,13\},\{4,11,12\},\{4,11,13\},\{11,12,13\}\}.$$
The configuration $\A$ and its regular triangulation $\Delta$ are
shown in Figure~\ref{ex1}. The toric ideal $I_{\A}$ lives in $R =
{\mathbb K}[a, \ldots, m]$. In Figure~\ref{ex1}, we have labeled the
points of $\A$ by the variables $a, \ldots, m$ corresponding to the
columns of ${\bf A}$, instead of by $1,\ldots, 13$.
\begin{figure}
\input{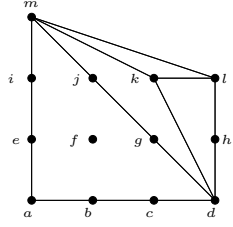}
\caption{The graded supernormal configuration $\A$ of
  Example~\ref{example1}.} \label{ex1} 
\end{figure}
The term order $\succ$ in Theorem~\ref{dnormal} can be induced 
via the weight vector $(7,5,3,1,5,3,1,1,3,1,0,1,1)$ refined by 
the reverse lexicographic order with $b > e > c > f > i > g > j > h > 
a > d > m > l > k$.  The initial ideal (computed using Macaulay~2
\cite{M2}) is $ J = \langle jl, gl, hm, h^2, j^2, gj, ik, fk, il, fl, jh,
cl, gh, ih, ch, ij, fj, ig, \\ ek, el, bl, fh, g^2, ck, bh,
cg, ej, i^2, fi, c^2, ak, f^2, ci, eg, al, eh, fg, cj, bk,
ha, cf, bg, ei, bi,\\ ef, bf, ec, bc, e^2, be, b^2, dml
\rangle \subset R.$ Its standard pairs,
grouped by the facets of $\Delta$ are:\\ 
$$\begin{array}{|c|c|}
\hline
\textup{faces} & \textup{roots}\\
\hline
\{1, 4, 13\} & 1, b, c, e, f, g, i, j, bj\\
\{ 4, 11, 13\} & 1,g,j\\
\{11, 12, 13\} & 1\\
\{4, 11, 12\} & 1, h\\
\hline
\end{array}$$\\
For $\sigma = \{1,4,13\}$, $FP_{\sigma}$ consists of nine lattice
points --- ${\bf Au}$ for each exponent vector ${\bf u}$ of the roots
$1, b, c, e, f, g, i, j, bj$. The last of these is $(2,2,2)^t$. The
monomials ${\bf x^u}$ of $R$ such that ${\bf Au} = (2,2,2)^t$ in
increasing order with respect to $\succ$ are: $bj,eg,ci,f^2, ak$.
Thus, $(bj,\{1,4,13\}) \in \mathcal S(J)$.





\end{example}

\section{Cohen-Macaulayness} \label{cmsection}

In this section we prove that the initial ideal $J$ of
Theorem~\ref{dnormal} is Cohen-Macaulay.  This is done by showing that
$J$ has a particular {\em Stanley filtration} \cite{ms} which implies
Cohen-Macaulayness \cite{simon,ms}. Stanley filtrations are special
{\em Stanley decompositions}.

\begin{definition}
  Let $M \subseteq R$ be a monomial ideal. A {\bf Stanley
    decomposition} of $M$ is a set of pairs of $M$, $\{ ({\bf x^u},
  \tau) \}$, that {\em partition} the standard monomials of $M$.
\end{definition}

\begin{remark}
  Every monomial ideal $M$ has the trivial Stanley decomposition $\{
  ({\bf x^u}, \emptyset): {\bf x^u} \notin M \}$. There can be many
  Stanley decompositions of a monomial ideal. The standard pair
  decomposition of $M$ is usually not a partition of the standard
  monomials of $M$.
\end{remark}

We will show that the standard pair decomposition $\mathcal S(J)$ of
$J$ can be modified first to a Stanley decomposition and then to a
Stanley filtration of the needed form. For $\tau \subseteq [n]$ let $
\pi_{\tau}: R \rightarrow R_{\tau} := {\mathbb K}[x_i: i \notin \tau]
$ be the projection map where $\pi_{\tau}(x_i) = x_i$ if $i \notin \tau$ and
$\pi_{\tau}(x_i) = 1$ if $i \in \tau$. The following is a well known
fact.

\begin{theorem} \cite[\S 12.D]{stu}\label{artproj}
  If $\sigma \in \max \, \Delta_{\omega}$ for a regular triangulation
  $\Delta_{\omega}$ of $\A$, then $\pi_{\sigma}(in_{\omega}(I_\A))$ is
  an artinian monomial ideal in $R_{\sigma}$ with $vol(\sigma)$-many
  standard monomials which are precisely the roots of the standard
  pairs $(\ast, \sigma)$ of $in_{\omega}(I_\A)$.
\end{theorem}

\begin{corollary} \label{rootsdivisors}
  If $({\bf x^u}, \sigma)$ is a standard pair in $\mathcal S(J)$, then every
  divisor of ${\bf x^u}$ is also the root of a standard pair in $\mathcal
  S(J)$. 
\end{corollary}

\begin{lemma} \label{rootslemma}
  Let $({\bf x^u}, \sigma)$ and $({\bf x^v}, \tau)$ be two standard pairs in
  $\mathcal S(J)$. If ${\bf x^u} \neq {\bf x^v}$ then $({\bf x^u}, \sigma) 
  \cap ({\bf x^v},\tau) = \emptyset$.
\end{lemma}

\begin{proof} Suppose ${\bf x^m} \in ({\bf x^u}, \sigma) \cap ({\bf x^v}, 
  \tau)$. 
  Then ${\bf x^m} = {\bf x^u} {\bf x}^{{\bf m}_\sigma}_\sigma = {\bf x^v} 
  {\bf x}^{{\bf m}_\tau}_\tau$ with
  $supp({\bf x}^{{\bf m}_\sigma}_\sigma) \subseteq \sigma$, 
  $supp({\bf x}^{{\bf m}_\tau}_\tau)
  \subseteq \tau$ and $supp({\bf x^u}), supp({\bf x^v})$ outside the vertices 
  of $\Delta_{\succ}$ and thus in particular, outside $\sigma \cup \tau$.
  Hence, ${\bf x^u} = {\bf x^v}$.
\end{proof}

\begin{corollary} \label{simpcase}
  If $\A$ is normal, $cone(\A)$ is simplicial (generated without loss
  of generality by ${\bf a}_1, \ldots, {\bf a}_d$), and $J$ is the special initial
  ideal of Theorem~\ref{dnormal} supported on $\Delta_{\succ} = \{\{1,
  \ldots, d\}\}$, then $\mathcal S(J)$ is a Stanley decomposition.
\end{corollary}

\begin{proof}
  Here $\A$ is $\Delta_{\succ}$-normal. All standard pairs in
  $\mathcal S(J)$ have face $[d]$ and roots the standard monomials of
  $\pi_{[d]}(J)$. Thus no two standard pairs of $\mathcal S(J)$ have
  the same root and by Lemma~\ref{rootslemma}, no two standard
  pairs intersect.
\end{proof}

However, if $|max \, \Delta_{\succ}| > 1$, then it is precisely the
standard pairs in $\mathcal S(J)$ with a common root that stop
$\mathcal S(J)$ from being a partition.  Such pairs always
exist when $|max \, \Delta_{\succ}| > 1$ --- for instance, $(1,
\sigma)$ is in $\mathcal S(J)$ for all $\sigma \in \textup{max} \,
\Delta_{\succ}$.  We will use the combinatorial notion of {\em
  shellings} to create new pairs that partition the standard
monomials covered by each set of standard pairs with a common root.
In Section~\ref{setup} we defined $ \sigma_{in} := \{ i \, : \,
{\bf a}_i \in cone(\A_{\sigma}), i \notin \sigma \}$. The following lemma 
identifies the faces in all standard pairs that share a root.

\begin{lemma} \label{sigmasets} 
  If ${\bf x^u}$ is a root of a standard pair in $\mathcal S(J)$, then $\{
  \sigma \in \textup{max} \, \Delta_{\succ}\, : \, ({\bf x^u}, \sigma) \in
  \mathcal S(J) \} = \{ \sigma \in \textup{max} \, \Delta_{\succ} \, :
  \, supp({\bf x^u}) \subseteq \sigma_{in} \}$.
\end{lemma}

\begin{proof}
  Recall from Section~\ref{setup} that if $({\bf x^u}, \sigma) \in \mathcal
  S(J)$, then $supp({\bf x^u}) \subseteq \sigma_{in}$. Conversely, suppose
  $({\bf x^u}, \tau)$ is a standard pair in $\mathcal S(J)$ and $supp({\bf x^u})
  \subseteq \sigma_{in}$ for some $\sigma \neq \tau$ in $\textup{max}
  \, \Delta_{\succ}$.  Then $supp({\bf x^u}) \subseteq \tau_{in} \cap
  \sigma_{in} = (\tau \cap \sigma)_{in}$. Then ${\bf Au} \in cone(\A_{\sigma
    \cap \tau})$. Since ${\bf Au} \in FP_{\tau} \cap \Z^d \cap
  cone(\A_{\sigma \cap \tau})$, it is also in $FP_{\sigma} \cap \Z^d$.
  Further, ${\bf x^u}$ is the cheapest monomial with respect to $\succ$
  among all monomials ${\bf x^v}$ in $R$ with ${\bf Au} = {\bf Av}$. This implies that
  $({\bf x^u}, \sigma)$ is also a standard pair of $J$.
\end{proof}

\begin{definition} \label{shelling} \cite[Chapter 3, Definition
  2.1]{cca} Let $\CC$ be a pure simplicial complex. A {\bf shelling}
  of $\CC$ is a linear ordering $F_1, \ldots, F_s$ of the facets of
  $\CC$ such that for each $j$, $1 < j \leq s$, the subcomplex
  supported in $(F_1 \cup \cdots \cup F_j) \backslash (F_1 \cup \cdots
  \cup F_{j-1})$ has a unique minimal face. A simplicial complex $\CC$
  is {\bf shellable} if it has a shelling.
\end{definition}

Let $F$ be a face of a simplicial complex $\CC$. Then $star(F,\CC) :=
\{ G \in \CC \, : \, F \cup G \in \CC \}$ is the simplicial complex
generated by all $G \in \CC$ containing $F$. We sometimes write
$star(F)$ for $star(F,\CC)$ when $\CC$ is obvious. The following is a
mild generalization of Lemma 8.7 in \cite{z}.

\begin{lemma} \label{starshell} Let $\CC$ be a pure shellable
  simplicial complex with shelling order $F_1, \ldots, F_s$. If $F$ is
  any face of $\CC$ then the restriction of the global shelling order
  to $star(F, \CC)$ yields a shelling of $star(F, \CC)$.
\end{lemma}

\begin{definition}
  For a root ${\bf x^u}$ in $\mathcal S(J)$, let $ \delta({\bf x^u}) := 
  \bigcap{\{ \sigma \, : \, ({\bf x^u}, \sigma) \in \mathcal S(J) \} }$.
\end{definition}

In the following arguments we fix a root ${\bf x^u}$ of a standard pair in
$\mathcal S(J)$.  By Lemma~\ref{sigmasets}, $ \delta({\bf x^u}) = \bigcap{
  \{ \sigma \in \textup{max} \, \Delta_{\succ} \, : \, supp({\bf x^u})
  \subseteq \sigma_{in} \} }$ and the set of facets 
$$ 
   \textup{max} \, star(\delta({\bf x^u}), \Delta_{\succ}) = 
   \{ \sigma \, : \, ({\bf x^u}, \sigma) \in \mathcal
   S(J) \} = \{ \sigma \in \textup{max} \, \Delta_{\succ} \, : \,
   supp({\bf x^u}) \subseteq \sigma_{in} \}. 
$$ 
If ${\bf x^v}$ divides ${\bf x^u}$, 
then $star( \delta({\bf x^u})) \subseteq star( \delta({\bf x^v}))$. The regular
triangulation $\Delta_{\succ}$ is shellable \cite{z}. In the rest of
this section, we fix a shelling of $\Delta_{\succ}$. By
Lemma~\ref{starshell}, this induces a shelling of $star(\delta({\bf x^u}))$.
Let us assume without loss of generality that $\sigma_{1},\sigma_{2},
\ldots, \sigma_{t}$ is the induced shelling order of the facets of
$\textup{max} \, star(\delta({\bf x^u}))$. For $\sigma_j \in \textup{max} \,
star(\delta({\bf x^u}))$, let $Q^{\sigma_j}_{\bf u}$ be the unique minimal face
described in Definition~\ref{shelling}. It is known that
$Q^{\sigma_j}_{\bf u} := \{ v \in \sigma_j \, : \,\sigma_j \backslash \{v\}
\, \subseteq \sigma_l \, \textup{for some}\, \, l < j \}.$ 
\vspace{.1cm}

\noindent{\bf Example~\ref{example1} continued.}
Consider the root $g$ of ${\mathcal S(J)}$ and the shelling order 
$\sigma_1 = \{4,11,12\}, \sigma_2 = \{11,12,13\}, \sigma_3
= \{4,11,13\}$ and $\sigma_4 = \{1,4,13\}$ on $\Delta_{\succ}$. 
Then $\sigma_3, \sigma_4$ is the induced shelling order of
$star(\delta(g))$. From this we obtain the sets $Q_g^{\sigma_3} = \emptyset$ 
and $Q_g^{\sigma_4} = \{1\}$.
\vspace{.1cm}

Consider the interval
$ I^{\sigma_j}_{\bf u} := \{ F \,:\, Q^{\sigma_j}_{\bf u} \subseteq F \subseteq
\sigma_j \}.$ 
\begin{lemma} \cite[pp. 247]{z} \label{localgeompart} 
  The simplicial complex $star( \delta ({\bf x^u}))$ is the disjoint union
  of the intervals $I^{\sigma_j}_{\bf u}$, $j = 1, \ldots, t$.
\end{lemma}

\begin{remark} \label{subshellingremark}
  Note that by construction, the partial union $I^{\sigma_1}_{\bf u} \cup
  I^{\sigma_2}_{\bf u} \cup \cdots \cup I^{\sigma_l}_{\bf u}$ is a partition
  of the subcomplex with maximal faces $\sigma_1, \sigma_2, \ldots,
  \sigma_l$ and that $Q^{\sigma_j}_{\bf u}$ is not contained in this
  partial union for any $j > l$.
\end{remark}

\begin{definition}
For a root ${\bf x^u}$ in $\mathcal S(J)$ and a facet $\sigma \in
\textup{max} \, \Delta_{\succ}$ define the monomial
$m^{\sigma}_{\bf u} :=  \left \{ \begin{array}{lll}
        1 &   \textup{if}\,\, \sigma \in star(\delta({\bf x^u})), \,
              Q_{\bf u}^{\sigma} = \emptyset \\
        \prod_{l \in Q^{\sigma}_{\bf u}}  x_l &  \textup{if}
        \,\, \sigma \in star(\delta({\bf x^u})), \, Q_{\bf u}^{\sigma} 
        \neq \emptyset \\
        1 & \textup{otherwise.}
                        \end{array} \right.$ 
\end{definition}

Recall that we have fixed a shelling of $\Delta_{\succ}$ and thus, 
by Lemma~\ref{starshell}, a
shelling of $star(\delta({\bf x^u}))$ for each root ${\bf x^u} \in \mathcal S(J)$.
Therefore, if $\sigma \in \textup{max} \, star(\delta({\bf x^u}))$,
$Q^{\sigma}_{u}$ is uniquely defined. We return to the fixed root 
${\bf x^u}$ and the shelling $\sigma_1, \ldots, \sigma_t$ of
$star(\delta({\bf x^u}))$. 

\begin{lemma} \label{partition}
  The standard monomials of $J$ in $\bigcup_{j = 1}^{t} {({\bf x^u},
    \sigma_j)}$ are partitioned by the pairs $({\bf x^u} \cdot
  m^{\sigma_j}_{\bf u}, \sigma_j)$, $j = 1, \ldots, t$. 
\end{lemma}

\begin{proof}
  By Lemma~\ref{localgeompart}, $I^{\sigma_1}_{\bf u} \cup
  I^{\sigma_2}_{\bf u} \cup \cdots \cup I^{\sigma_t}_{\bf u}$ is a partition
  of $star( \delta ( {\bf x^u}))$.  Hence, if ${\bf x^v} \in \bigcup_{j = 1}^{t} {
    ({\bf x^u}, \sigma_j ) }$ then $supp({\bf x^{v-u}} ) \in I^{\sigma_j}_{\bf u}$ for
  a unique $I^{\sigma_j}_{\bf u}$. By construction, $I^{\sigma_j}_{\bf u} = \{
  F \, : \, supp(m^{\sigma_j}_{\bf u}) \subseteq F \subseteq \sigma_j \}$
  and so ${\bf x^v} = {\bf x^u} \cdot {\bf x^{v-u}} \in ({\bf x^u} \cdot 
  m^{\sigma_j}_{\bf u},
  \sigma_j)$. This implies that $\bigcup_{j = 1}^{t} { ({\bf x^u},
    \sigma_j ) } = \bigcup_{j = 1}^{t} { ({\bf x^u} \cdot
    m^{\sigma_j}_{\bf u} , \sigma_j ) } $ where the inclusion $\subseteq$
  follows from the previous line and $\supseteq$ from the fact that
  $({\bf x^u} \cdot m^{\sigma_j}_{\bf u} , \sigma_j) \subseteq ({\bf x^u},
  \sigma_j)$ for each $j = 1, \ldots, t$. To see that $\bigcup_{j =
    1}^{t} { ({\bf x^u} \cdot m^{\sigma_j}_{\bf u} , \sigma_j ) }$ is a
  partition, suppose ${\bf x^v} \in ({\bf x^u} \cdot m^{\sigma_i}_{\bf u} ,
  \sigma_i ) \cap ({\bf x^u} \cdot m^{\sigma_j}_{\bf u} , \sigma_j )$ where
  $i < j$. Then ${\bf x^{v-u}}$ has support in $I^{\sigma_i}_{\bf u} \cap
  I^{\sigma_j}_{\bf u} = \emptyset$ which implies that ${\bf x^v} = {\bf x^u}$. However,
  for $j > 1$, $m^{\sigma_j}_{\bf u} \neq 1$ as $Q^{\sigma_j}_{\bf u} \neq
  \emptyset$ which means that ${\bf x^u}$ lies only in $({\bf x^u}, \sigma_1)$.
\end{proof}

\noindent{\bf Example~\ref{example1} continued.}
  As before, the monomial $g$ is a root of ${\mathcal S(J)}$ with the shelling 
  order induced on $star(\delta (g) )$ as above. We obtain the monomials 
  $m^{\sigma_3}_g = 1$ and  $m^{\sigma_4}_g = a$ from 
  $Q^{\sigma_3}_g = \emptyset$ and  $Q^{\sigma_4}_g = \{ 1 \} $. 
  Then $(g, \{4,11,13\}) \cup (g,\{1,4,13\}) 
  = (g, \{4,11,13\}) \cup (g \cdot a,\{1,4,13\})$ with the latter union 
  being disjoint.
\begin{theorem} \label{stanley}
  Let $\sigma_1, \ldots, \sigma_s$ be the fixed shelling of
  $\Delta_{\succ}$. Then
\begin{equation} 
\bigcup_{i=1}^{s} \bigcup_{({\bf x^u},\sigma_i) \in \mathcal S(J)} 
({\bf x^u} \cdot m_{\bf u}^{\sigma_i}, \sigma_i) 
\label{standecomp}
\end{equation}
is a Stanley decomposition of $J$.
\end{theorem}

\begin{proof}
  Lemma~\ref{partition} showed how to make the union of the standard
  pairs of $\mathcal S(J)$ with a common root a disjoint union of
  pairs of $J$. By Lemma~\ref{rootslemma}, (\ref{standecomp}), the
  union of these disjoint unions is a Stanley decomposition of
  $J$.
\end{proof}

The above Stanley decomposition can be organized to have more
structure.

\begin{definition} \label{filtration}~\cite{ms}
  Let $M \subseteq R$ be a monomial ideal. A {\bf Stanley filtration}
  of $M$ is a Stanley decomposition of $M$ with an ordering of the
  pairs $ \{ ({\bf x}^{ {\bf v}_i}, \tau_i) \, : \, 1 \leq i \leq r \} $ such
  that for all $1 \leq j \leq r$ the set $ \{ ({\bf x}^{ {\bf v}_i}, \tau_i)
  \, : \, 1 \leq i \leq j \} $ is a Stanley decomposition of $M_j :=
  M + \langle {\bf x}^{ {\bf v}_{j+1}}, {\bf x}^{ {\bf v}_{j+2}}, \ldots , {\bf x}^{
    {\bf v}_r} \rangle$.  Equivalently, the ordered set $ \{ ({\bf x}^{{\bf v}_i},
  \tau_i) \, : \, 1 \leq i \leq r \} $ is a Stanley filtration
  provided the modules $R/{M_j}$ form a filtration $ {\mathbb K} = R/{M_0}
  \subsetneq R/{M_1} \subsetneq R/{M_2} \subsetneq \cdots \subsetneq
  R/{M_r} = R/M $ with $\frac{R/{M_j}}{R/{M_{j-1}}} \cong {\mathbb K}
  [x_i: i  \in \tau_j]$.
\end{definition}

\begin{example}
  (from~\cite{ms}) Let $M = \langle x_1x_2x_3 \rangle \subset {\mathbb
    K}[x_1, x_2, x_3]$. Then $$\{ (1, \emptyset), (x_1, \{1,2\}),
  (x_2, \{2,3\}), (x_3, \{1,3\}) \}$$
  is a Stanley decomposition of 
  $M$ but no ordering of these pairs is a Stanley filtration of $M$.
  Alternatively, the ordered pairs $(1, \{ 1,3 \} ), (x_2,
  \{2,3\}), (x_1x_2, \{1,2\})$ form a Stanley filtration of $M$.
\end{example}

We now show that the pairs in (\ref{standecomp}) can be ordered to
yield a Stanley filtration of $J$. The significance of this for us
comes from a result of Simon~\cite{simon}, interpreted as follows by
Maclagan and Smith~\cite{ms}.

\begin{theorem} \label{filtgivescm} If $M \subseteq R$ is a monomial
  ideal with a Stanley filtration such that for each face $\tau$ of a
  pair in the filtration, the prime ideal $P_{\tau}$ is a minimal
  prime of $M$, then $M$ is Cohen-Macaulay.
\end{theorem}

The faces of pairs in (\ref{standecomp}) already index minimal primes
of $J$. Thus to show that $J$ is Cohen-Macaulay all we need to do is
to order the pairs in (\ref{standecomp}) so that the ordered
decomposition is a Stanley filtration. We do this using the following
algorithm.

\begin{algorithm} \label{list}
{\bf Input:} The Stanley decomposition (\ref{standecomp}) of $J$.\\
{\bf Output:} A Stanley filtration of $J$ with the same faces as those in
(\ref{standecomp}). 
\begin{description}
\item[1] (Local Lists) For each $\sigma_i$, $1 \leq i \leq s$, 
order the pairs in (\ref{standecomp}) with face $\sigma_i$
in any way such that if $({\bf x^u} \cdot m_{\bf u}^{\sigma_i}, \sigma_i)$ 
precedes $({\bf x^v} \cdot m_{\bf v}^{\sigma_i}, \sigma_i)$ then ${\bf x^v}$ does not 
divide ${\bf x^u}$. Call this list $L_i$.
\item[2] (Global List) The global list $\LL$ is obtained by 
appending $L_i$ to the end of $L_{i-1}$ for $i = 2, \ldots, s$.
\end{description}
\end{algorithm}

\begin{proof} Let $r_i := \sum_{l=1}^i
  {vol(\sigma_l)}$ for $i = 1, \ldots, s$. Then $r := r_s$ is the
  total number of pairs in (\ref{standecomp}).  Write $\LL$ as
  $[({\bf x}^{{\bf u}_l} \cdot m^{\tau_l}_{{\bf u}_l} , \tau_l) \,:\, 1 \leq l \leq
  r ]$ where $\tau_l = \sigma_i$ when $r_{i-1} < l \leq r_i$ ($r_0 :=
  0$) and ${\bf x}^{{\bf u}_l} \cdot m^{\tau_l}_{{\bf u}_l}$ is the root of the
  $(l-{r_{i-1}})$-th pair in the local list $L_i$ constructed in
  Step 1 of the algorithm. For $1 \leq j \leq r$ define the partial list $\LL_j :=
  [({\bf x}^{{\bf u}_l}\cdot m^{\tau_l}_{{\bf u}_l} , \tau_l) \,:\, 1 \leq l \leq
  j]$ and the ideal $M_j := J + \langle \, {\bf x}^{ {\bf u}_{j+1}} \cdot
  m^{\tau_{j+1}}_{{\bf u}_{j+1}} \, , \, {\bf x}^{ {\bf u}_{j+2}} \cdot
  m^{\tau_{j+2}}_{{\bf u}_{j+2}} \, , \, \ldots , \, {\bf x}^{ {\bf u}_r} \cdot
  m^{\tau_r}_{{\bf u}_{r}} \, \rangle$. We need to prove that $\LL_j$ is a
  Stanley decomposition of $M_j$. Since $\LL_j$ is already a
  partition, it suffices to show that the set of monomials in the
  pairs in $\LL_j$ is the set of standard monomials of $M_j$.
  
  (i) {\em The standard monomials of $M_j$ are contained in the pairs
    in $\LL_j$}: A standard monomial ${\bf x^u}$ of $M_j$ is a standard
  monomial of $J$ and hence is covered by a unique pair $({\bf x}^{ {\bf u}_l}
  \cdot m^{\tau_l}_{{\bf u}_l} , \tau_l)$ in $\LL$. Also, ${\bf x^u} \notin
  \langle {\bf x}^{ {\bf u}_{j+1}} \cdot
  m^{\tau_{j+1}}_{{\bf u}_{j+1}},\ldots,{\bf x}^{{\bf u}_r} \cdot
  m^{\tau_r}_{{\bf u}_{r}} \rangle$ which implies that ${\bf x^u} \notin
  ({\bf x}^{{\bf u}_{j+k}} \cdot m_{{\bf u}_{j+k}}^{\tau_{j+k}}, \tau_{j+k})$ for any $k
  \geq 1$. Hence $l \leq j$ and ${\bf x^u} \in ({\bf x}^{{\bf u}_l} \cdot
    m^{\tau_l}_{{\bf u}_l},\tau_l) \in \LL_j$.
  
    (ii) {\em The monomials in the pairs in $\LL_j$ are standard
      monomials of $M_j$}: Suppose ${\bf x^u}$ lies in the (unique) pair
    $({\bf x}^{{\bf u}_l} \cdot m_{{\bf u}_l}^{\tau_l}, \tau_l) \in \LL_j$.  
    Since ${\bf x^u}
    \notin J$, it suffices to show that ${\bf x^u} \notin \langle
    {\bf x}^{{\bf u}_{j+1}} \cdot m^{\tau_{j+1}}_{{\bf u}_{j+1}}, \ldots, 
    {\bf x}^{{\bf u}_r} \cdot m^{\tau_r}_{{\bf u}_r}\rangle$.
   
    Suppose ${\bf x^u} \in \langle {\bf x}^{{\bf u}_{j+1}} \cdot
    m^{\tau_{j+1}}_{{\bf u}_{j+1}}, \ldots, {\bf x}^{{\bf u}_r} 
    \cdot m_{{\bf u}_r}^{\tau_r} \rangle$. 
    Then there exists $p$, $j+1 \leq p \leq r$ such that
    ${\bf x}^{{\bf u}_p} \cdot m_{{\bf u}_p}^{\tau_p} \, | \, {\bf x^u} = 
    {\bf x}^{{\bf u}_l}\cdot m_{{\bf u}_l}^{\tau_l} \cdot {\bf x}_{\tau_l}^{\ast}$ where
    ${\bf x}_{\tau_l}^{\ast}$ is a monomial with support in $\tau_l$.  Since
    $supp({\bf x}^{{\bf u}_p})$ and $supp({\bf x}^{{\bf u}_l})$ are both in 
    $[n] \backslash (\tau_p \cup \tau_l)$, it follows that 
    ${\bf x}^{{\bf u}_p} | {\bf x}^{{\bf u}_l}$.  Since $l < p$, by Step 1 of 
    the algorithm, $\tau_p \neq \tau_l$.  Recall
    that $({\bf x}^{{\bf u}_l}, \tau_l)$ and $({\bf x}^{{\bf u}_p}, \tau_p)$ are standard
    pairs of $J$. Since ${\bf x}^{{\bf u}_p} | {\bf x}^{{\bf u}_l}$, by
    Corollary~\ref{rootsdivisors}, $({\bf x}^{{\bf u}_p}, \tau_l)$ is also in
    $\mathcal S(J)$. This implies that $\tau_p$ and $\tau_l$ are two
    distinct facets in $star(\delta({\bf x}^{{\bf u}_p}))$. Since
    $m_{{\bf u}_p}^{\tau_p} | {\bf x^u}$, $Q_{{\bf u}_p}^{\tau_p} (=
    supp(m_{{\bf u}_p}^{\tau_p})) \subseteq supp({\bf x^u}) \cap \bigcup_{i=1}^{s}
    \sigma_i \subseteq \tau_l$. However, this is a contradiction since
    $\tau_l$ precedes $\tau_p$ in the shelling order on
    $\Delta_{\succ}$ and hence $Q_{{\bf u}_p}^{\tau_p}$ cannot be in
    $\tau_l$. Thus $m_{{\bf u}_p}^{\tau_p} \not | {\bf x^u}$ and ${\bf x^u} \notin 
    \langle {\bf x}^{{\bf u}_{j+1}} \cdot
    m^{\tau_{j+1}}_{{\bf u}_{j+1}}, \ldots, {x}^{{\bf u}_r} \cdot m_{{\bf u}_r}^{\tau_r}
    \rangle$ and thus not in $M_j$.
\end{proof}

\noindent{\bf Example~\ref{example1} continued.}
As before, $\sigma_1 = \{4,11,12\}, \sigma_2 = \{11,12,13\}, \sigma_3
= \{4,11,13\}$ and $\sigma_4 = \{1,4,13\}$ is a shelling order on
$\Delta_{\succ}$. The (ordered) local lists in the Stanley filtration
produced by Algorithm~\ref{list} are: \\
$\mathcal L_1 = [(1,\{4,11,12\}), (h,\{4,11,12\})]$,\\
$\mathcal L_2 = [(1 \cdot m, \{11,12,13\})]$,\\
$\mathcal L_3 = [(1 \cdot dm, \{4,11,13\}), (g,\{4,11,13\}), (j,
\{4,11,13\})]$, \\
$\mathcal L_4 = [(1 \cdot a, \{1,4,13\}), (b,\{1,4,13\}), (c,
\{1,4,13\}), (e, \{1,4,13\}), (f, \{1,4,13\}),\\ (g \cdot a,
\{1,4,13\}), (i, \{1,4,13\}), (j \cdot a, \{1,4,13\}),
(bj,\{1,4,13\})]$.
\begin{theorem} \label{cm}
  Let $\A$ be a $\Delta$-normal configuration. Then there exists a
  term order $\succ$ such that $\Delta = \Delta_{\succ}$ and
  $in_{\succ}(I_\A)$ is Cohen-Macaulay.
\end{theorem}

\begin{proof} Algorithm~\ref{list} shows that the initial ideal $J$ of
  Theorem~\ref{dnormal} has a Stanley filtration that satisfies the
  conditions of Theorem~\ref{filtgivescm}. This theorem guarantees
  that $J$ is Cohen-Macaulay.
\end{proof}

\begin{remark}
  We remark that even when $\A$ is $\Delta$-normal it is not true that
  all initial ideals of $I_{\A}$ without embedded primes are
  Cohen-Macaulay. Take $\A$ to be the columns of
  $${\bf A} = \left ( \begin{array}{cccccccc}
      1 & 1 & 1 & 1 & 1 & 1 & 1 & 1\\
      1 & 0 & 7 & 4 & 3 & 2 & 5 & 4\\
      0 & 1 & 5 & 7 & 5 & 4 & 5 & 4\\
      0 & 0 & 6 & 6 & 4 & 3 & 5 & 4 \end{array} \right ).$$
  Then $\A$
  admits a unimodular regular triangulation and is hence
  $\Delta$-normal. The toric ideal $I_{\A} \subset {\mathbb
    K}[a,\ldots,h]$ has codimension four and has 46 initial ideals
  without embedded primes. Among them, the following two have
  projective dimension five. 
\begin{enumerate}
\item $ \langle acd,adg,afg,ae,ag^2,ce,cf,eh,f^2,bc^2d,fgh
  \rangle$
\item $ \langle acd,adg,afg,ae,ag^2,ce,cf,eh,f^2,fgh,g^2h^2
 \rangle $
\end{enumerate}
The initial ideals of $I_\A$ were computed using the software 
package CaTS \cite{cats} and then checked for embedded primes and 
Cohen-Macaulayness using Macaulay~2. 

We remark that the first example of a monomial toric initial ideal
without embedded primes that is not Cohen-Macaulay was found by Laura
Matusevich~\cite{lm}.  In that example, $I_{\A}$ is not Cohen-Macaulay and thus
$\A$ is not normal.

\end{remark}

\section{Degree Bounds} \label{dnormalboundsection}

\begin{theorem} \label{dnormalbound}
  If $\A$ is a graded $\Delta$-normal configuration, then there exists
  a term order $\succ$ such that $\Delta = \Delta_\succ$ and the
  Gr\"obner basis of $I_{\A}$ with respect to $\succ$ consists of
  binomials of degree at most $d = dim(\A)$.
\end{theorem}

Theorem~\ref{dnormalbound} settles Conjecture~\ref{bigconjecture} for
the subset of normal configurations that are $\Delta$-normal. Since
$\A$ is graded, $I_{\A}$ is homogeneous with respect to the usual
grading of $R$ where $deg(x_i) = 1$ for $i = 1, \ldots, n$. Hence
it suffices to show that $I_{\A}$ has an initial ideal of degree at
most $d$. We will show that the initial ideal $J$ from
Theorem~\ref{dnormal} satisfies this degree bound when $\A$ is graded.
Proposition~13.15 in \cite{stu} shows that
Conjecture~\ref{bigconjecture} is true whenever $\A$ admits a regular
unimodular triangulation. (See also Proposition~13.18 in \cite{stu}.)
Such configurations form a proper subset of the set of $\Delta$-normal
configurations.  If we are allowed to replace a graded normal $\A$ by
all the lattice points in a ``big enough'' multiple of the convex hull
of $\A$, then it is known that this new configuration admits regular
unimodular triangulations and thus has a Gr\"obner basis of degree at
most $d$.  See \cite{bgt1} and \cite{bgt2} for many such results.

Conjecture~\ref{bigconjecture} requires that $\A$ be both graded and
normal.
\begin{example} \label{needgradedandnormal}
  {\em Graded, but not normal:} When $\A = \{(1,0), (1,p), (1,q)\}$
  with $0 < p < q$, $q > 2$ and $g.c.d(p,q) = 1$, then $I_{\A} =
  \langle x_1^{q-p} x_3^p - x_2^q \rangle$. Its two initial ideals
  are therefore generated in degree $q > 2 = d$.\\
  {\em Normal, but not graded:} The normal configuration $\A = \{
  (1,0), (1,1), (p,p+1) \}$ where $p \geq 2$ has the toric ideal
  $I_{\A} = \langle x_1x_3 - x_2^{p+1} \rangle$. Hence $x_1x_3 -
  x_2^{p+1}$ is the unique element in both its reduced Gr\"obner
  bases.
\end{example}

\begin{remark} \label{sharp}~ (\cite{stu}, Chapter 13)
  The bound in Conjecture~\ref{bigconjecture} is best possible.
  Consider the graded $\Delta$-normal configuration $\A = \{ d {\bf
    e}_1, \, d {\bf e}_2, \ldots, d {\bf e}_d, \, {\bf e}_1 + {\bf
    e}_2 + \cdots + {\bf e}_d \}$ where $d \in \N$.  (Note that
  $cone(\A)$ simplicial).  Then $I_\A = \langle x_1 x_2 \cdots x_d -
  x_{d+1}^d \rangle$. 
\end{remark}

Consider the initial ideal $J$ from Theorem~\ref{dnormal} for a graded
$\Delta_{\succ}$-normal $\A$.  Since $\A$ is graded, we may assume
without loss of generality that ${\bf a}_i = (1, {\bf a}_i') \in \Z^d$
for $i = 1, \ldots, n$. We will show that $J$ is generated in degree
at most $d$.

For a $\sigma \in \textup{max} \, \Delta_{\succ}$, recall that
$\sigma_{in} := \{i \,: \, {\bf a}_i \in cone(\A_{\sigma}), \, i
\notin \sigma \, \}$. Define $\sigma_{out}:= \{i \,: \, {\bf a}_i
\notin cone(\A_{\sigma}) \}$. Then $\sigma \cup \sigma_{in} \cup
\sigma_{out}$ is a partition of $[n]$.  Let $J^{\sigma} :=
\pi_{\sigma}(J)$ be the artinian ideal in $ R_{\sigma} = {\mathbb
  K}[x_j \, : \, j \in \sigma_{in} \cup \sigma_{out}]$ from
Theorem~\ref{artproj}. Recall that the standard monomials of
$J^{\sigma}$ are the roots of standard pairs in $\mathcal S(J)$ with
face $\sigma$. Since the supports of these roots lie in $\sigma_{in}$,
$J^{\sigma} \cap {\mathbb K}[x_i \, : \, i \in \sigma_{in}]$ is a
monomial ideal $N^{\sigma} = \langle {\bf x}^{{\bf v}_1},{\bf x}^{
  {\bf v}_2},\ldots, {\bf x}^{{\bf v}_{r_\sigma}} \rangle$ with
$supp({\bf x}^{{\bf v}_i}) \subseteq \sigma_{in}$, and
$$J^{\sigma} = \langle x_j : j \in \sigma_{out}
\rangle \, + \, N^{\sigma}.$$

\begin{lemma} \label{lemma1ofdnormalbound}
  Each minimal generator ${\bf x}^{{\bf v}_i}$ of $N^{\sigma}$ is a
  minimal generator of $J$ of degree at most~$d$.
\end{lemma}
\begin{proof} A minimal generator ${\bf x}^{{\bf v}_i}$ of
  $N^{\sigma}$ is the projection via $\pi_{\sigma}$ of a minimal
  generator ${\bf x}^{{\bf v}_i} {{\bf x}^{\bf m}_\sigma}$ of $J$
  where $supp({{\bf x}^{\bf m}_\sigma}) \subseteq \sigma$. Suppose
  $supp({{\bf x}^{\bf m}_\sigma}) \not= \emptyset$. Then ${\bf x}^{{\bf v}_i}$ 
  is a standard monomial of $J$ with $supp({\bf x}^{{\bf v}_i}) \subseteq \sigma_{in}$. 
  Hence ${\bf x}^{{\bf v}_i}$ is
  covered by a standard pair $({\bf x}^{{\bf u}_\gamma}, \sigma)$ of
  $J$.  This implies that all monomials of the form ${\bf x}^{{\bf v}_i} 
  {\bf x}_{\sigma}^{\bf p}$ as ${\bf p}$ varies are standard
  monomials of $J$ which contradicts that ${\bf x}^{{\bf v}_i} { {\bf
    x}^{\bf m}_\sigma}$ is in $J$. Thus $supp({{\bf x}^{\bf
    m}_\sigma}) = \emptyset$ which implies that ${{\bf x}^{{\bf
      v}_i}}$ is a minimal generator of~$J$.
  
  Since ${\bf a}_i = (1, {\bf a}_i') \in \Z^d$ for $i \in [n]$, each
  lattice point in the half open fundamental parallelopiped
  $FP_{\sigma}$ of $cone(\A_{\sigma})$ lies on one of the $d$
  hyperplanes $x_1=0, \ldots, x_1 = d-1$ in $\R^d$.  Therefore, if
  ${\bf \gamma} \in FP_\sigma \cap \Z^d$, then the $1$-norm of ${\bf
    {u}_\gamma}$ which equals the first co-ordinate of $({\bf
    Au_{\gamma}})$ which equals ${\bf \gamma}_1$ is at most $d-1$.
  This implies that $deg({\bf x}^{\bf u_{\gamma}}) \leq d-1$. Thus all
  standard monomials of the artinian ideal $J^{\sigma}$ have degree at
  most $d-1$ which implies that the minimal generators of $J^{\sigma}$
  (and $N^{\sigma}$) have degree at most $d$.
\end{proof}

\noindent{\bf Example~\ref{example1} continued.}
For $\sigma = \{1,4,13\}$, $J^{\sigma} = \langle h, k, l \rangle +
(N^{\sigma} = \langle j^2, gj, ij, fj, ig,\\ g^2, cg, ej, i^2, fi, c^2,
f^2, ci, eg, fg, cj, cf, bg, ei, bi, ef, bf, ec, bc, e^2, be, b^2
\rangle)$. Note that all minimal generators of $N^{\sigma}$ are 
minimal generators of $J$ of degree at most three.

\begin{theorem} \label{simplexbound}
  If $\A$ is a graded normal configuration with $cone(\A)$ simplicial
  then $I_\A$ has a Gr\"obner basis consisting of binomials of degree
  at most $d$.
\end{theorem}

\begin{proof} 
  Assuming that $cone(\A)$ is generated by ${\bf a}_1, \ldots, {\bf
    a}_d$, $\A$ is $\Delta_{\succ}$-normal where $\Delta_{\succ}$ is
  the regular triangulation of $\A$ with the unique facet $\sigma =
  [d]$.  Here $\sigma_{out} = \emptyset$.
  
  We argue that all minimal generators of $J$ have support in
  $\sigma_{in} = \{d+1, \ldots, n\}$.  Suppose ${\bf x^{\alpha}}$ is a
  minimal generator of $J$ with $supp({\bf x^{\alpha}}) \cap [d] = F
  \neq \emptyset$. Let $G = supp({\bf x^{\alpha}}) \backslash [d]$.
  Then $G \neq \emptyset$ since otherwise ${\bf x^{\alpha}}$ would lie
  on the standard pair $(1, [d])$ of $J$ which is a contradiction.
  Write ${\bf x^{\alpha}} = {\bf x}^{{\bf \alpha}_F} {\bf x}^{{\bf
      \alpha}_G}$ where $supp({\bf \alpha}_F) \subseteq F$ and
  $supp({\bf \alpha}_G) \subseteq G$. Since $G, F \neq \emptyset$,
  ${\bf x}^{{\bf \alpha}_G}$ is a standard monomial of $J$ which
  implies that ${\bf x}^{\bf \alpha}$ is also a standard monomial of
  $J$ as ${\bf x}^{{\bf \alpha}_G}$ lies on some standard pair with
  face $[d]$.  This is a contradiction and so $F = \emptyset$.
  
  The above argument shows that $J$ and $N^{\sigma}$ have the same
  minimal generators. The degree bound then follows from the proof of
  Lemma~\ref{lemma1ofdnormalbound}.
\end{proof}

Theorem~\ref{simplexbound} proves Theorem~\ref{dnormalbound} in the
case where $cone(\A)$ is simplicial. When $cone(\A)$ is not
simplicial, $J$ may have minimal generators that are not pre-images
under $\pi_{\sigma}$ of the minimal generators of $N^{\sigma}$ (or
even $J^{\sigma}$) as $\sigma$ varies in $\textup{max} \,
\Delta_{\succ}$. Our next step is to show that for a $\sigma \in
\textup{max} \, \Delta_{\succ}$, the minimal generators of $J$ that
project under $\pi_{\sigma}$ to the minimal generators $x_j \in
\sigma_{out}$ of $J^{\sigma}$ have degree at most $d$. We need a
preliminary lemma.

Let $Q$ be a $(d-1)$-polytope in $ \{{\bf x} \in \R^d : x_1 = 0 \}$
and let $C$ be the cone over $Q$. Then there exists a matrix ${\bf S}
\in \R^{f \times d}$ such that $C = \{ {\bf x} \in \R^d \, : \, {\bf
  Sx} \geq 0 \}$ where each row of ${\bf S}$ is the normal to a facet
of $C$.  Hence $Q = \{ {\bf x} \in \R^d \, : \, x_1 = 1, \, {\bf Sx}
\geq {\bf 0} \}$.  Let $Q_{rev}$ be the system obtained by reversing
all the inequalities in $Q$:
$$
Q_{rev} = \{ {\bf x} \in \R^d \, : \, x_1 \leq 1, \, -x_1 \leq -1,
\, {\bf Sx} \leq {\bf 0} \}.
  $$
 
\begin{lemma} \label{halfspaces}
The polyhedron defined by $Q_{rev}$ is empty.
\end{lemma}

\begin{proof} We may assume that $Q$ has been translated so that the
  unit vector ${\bf e}_1 \in \R^d$ lies in the relative interior of
  $Q$.  If ${\bf x} \in C$ then by our assumption, $x_1 \geq 0$ which
  implies that ${\bf e}_1 \cdot {\bf x} (= x_1) \geq 0$.  This implies
  that ${\bf e}_1 \in C^* = \{{\bf yS} \, : \, {\bf y} \geq {\bf 0} \}$
where $C^*$ is the dual cone to $C$.  (Recall $C^* := \{{\bf v} \in
\R^d \, : \, {\bf v} \cdot {\bf x} \geq 0, \, \textup{for all} \, {\bf
  x} \in C \}$.)  Thus there exists some ${\bf y} \geq {\bf 0}, {\bf
  y} \neq {\bf 0}$ such that ${\bf yS} = {\bf e}_1$.  Therefore, if we
choose ${\bf v} \in \R^{2+f}$ such that ${\bf v} = (0,1,{\bf y})$ then
  ${\bf v} \geq {\bf 0}$, ${\bf v} \neq {\bf 0}$ and
  $$  {\bf v} \cdot \left(
\begin{array}{cccc} 1 & 0 & \cdots & 0 \\ -1 & 0 & \cdots & 0 \\
s_{11} & s_{12} & \cdots & s_{1d} \\ \vdots & \vdots & \vdots &
\vdots \\ s_{f1} & s_{f2} & \cdots & s_{fd} \end{array} \right ) = 0.$$ 
Let $z = (1,-1,0,\ldots,0)$ be the right hand side vector in 
the description of $Q_{rev}$. Then ${\bf v} \cdot {\bf z} = 1(-1) = -1
< 0$ and by Farkas lemma \cite[Prop. 1.7.]{z}, $Q_{rev} = \emptyset$.
\end{proof}

\begin{lemma} \label{lemma2ofdnormalbound}
  Let $\sigma$ be a facet of $\Delta_{\succ}$. Then for a $j \in
  {\sigma}_{out}$, the minimal generators of $J$ that are preimages of
  the minimal generator $x_j$ of $J^{\sigma}$ under the map
  $\pi_{\sigma}$ are squarefree monomials of degree at most~$d$.
\end{lemma}

\begin{proof}
  Let $\sigma \in \textup{max} \, \Delta_{\succ}$, $j \in
  {\sigma}_{out}$ and $P := x_j {\bf x}_{\sigma}^{\bf m}$ be a
minimal generator of $J$ with $Y := supp({\bf x}_{\sigma}^{\bf m})
\subseteq \sigma$.  All minimal generators of $J$ that project to
$x_j$ under $\pi_{\sigma}$ look like $P$.  If $Y = \emptyset$, then
$x_j$ is the only minimal generator of $J$ that projects to $x_j$ and
we are done.  Therefore, we consider the case where $Y \neq
\emptyset$.
  
Suppose $P$ is not squarefree. Then there exists an $i \in \sigma$
such that $m_i > 1$ where $m_i$ is the $i$-th co-ordinate of ${\bf
  m}$.  Since $P$ is a minimal generator of $J$, $P/x_i$ is a standard
monomial of $J$ with $supp(P/x_i) = supp(P) = \{j\} \cup Y$. Hence
there exists $\tau \in \textup{max} \, \Delta_{\succ}$ such that a
standard pair with face $\tau$ covers $P/x_i$. This implies that
$supp(P/x_i) = \{j\} \cup Y \subseteq \tau_{in} \cup \tau$. Since $Y
\subseteq \sigma$, $Y \cap \tau_{in} = \emptyset$ and thus, $Y
\subseteq \tau \cap \sigma$. If $j \in \tau$, then $P$ is covered by
the standard pair $(1, \tau)$ which contradicts that $P$ is in $J$.
If $j \in \tau_{in}$, then ${\bf a}_j$ lies in $cone(\A_{\tau})$.
Since $\A$ lies on the hyperplane $\{{\bf x} \in \R^d \, : \, x_1 = 1
\}$, ${\bf a}_j$ is in fact in the minimal Hilbert basis of both
$cone(\A_\tau)$ and $cone(\A)$ and hence ${\bf e}_j$ is the unique
vector in $\N^n$ that satisfies ${\bf Ax} = {\bf a}_j$.  Consequently
$(x_j, \tau)$ is a standard pair of $J$.  But this implies that $P$
lies on this standard pair which is again a contradiction.  Therefore,
$P$ is squarefree.
  
  To argue that $deg(P) \leq d$, it therefore suffices to prove
  that $Y \subsetneq \sigma$. Suppose $\sigma = [d]$, ${\bf x}_{[d]}
  := \prod_{i \in \sigma} x_i$ and $P = x_j {\bf x}_{[d]}$. Then for
  each $i \in [d]$, $P/x_i$ is a standard monomial of $J$ and is
  therefore covered by a standard pair $(\ast, \tau^i)$ of $J$. The
  face $\tau^i$ does not contain $i$ since otherwise $P$ would be a
  standard monomial of $J$. In particular, $\tau^i \neq [d]$ for any
  $i \in [d]$. Also, $j \in \tau^i \cup \tau^i_{in}$ for each $i \in
  [d]$.
  
  We now show that we may assume $\tau^i \cap [d] = [d]\backslash \{ i
  \}$ for all $i \in [d]$. Clearly, $\tau^i \cap [d] \subseteq [d]
  \backslash \{ i \}$ since $i \notin \tau^i$. Suppose a monomial in
  $P/x_i \cdot {\mathbb K}[x_l : l \in [d] \backslash \{ i \}]$ lies
  in $J$. Then it is divisible by a minimal generator of $J$ that
  projects to $x_j$ under $\pi_{\sigma}$, all of which are squarefree.
  Such a minimal generator would properly divide $P$ which contradicts
  that $P$ is a minimal generator of $J$. Hence $(P/x_i, [d]
  \backslash \{i\})$ is a pair of $J$ and therefore, contained in a
  standard pair of $J$. We may assume that $\tau^i$ is the face of
  this standard pair. Thus $[d] \backslash \{i\} \subset \tau^i$ and
  $\tau^i \cap [d] = [d]\backslash \{ i \}$ as claimed.
  
  Since $\A$ is graded, $\tau^1, \ldots, \tau^d$ index
  $(d-1)$-simplices in a regular triangulation of $conv(\A)$, the
  convex hull of $\A$. The simplex indexed by $[d]$ is geometrically
  $Q[d] = \{{\bf x} \in \R^d \,: \, {\bf s}_i \cdot {\bf x} \geq 0, \,
  x_1 = 1 , \, i = 1,2,\ldots, d \} $ where ${\bf s}_i \cdot {\bf a}_l
  = 0$ for all $l \in [d] \backslash \{i\}$ and ${\bf s}_i \cdot {\bf
    a}_i > 0$.  Now $j \in \tau^i \cup \tau^i_{in} \cap \sigma_{out}$
  for each $i \in [d]$ implies that ${\bf a}_j \in Q[d]_{rev}$ where $
  Q[d]_{rev} = \{ {\bf x} \in \R^d \,: \, {\bf s}_i \cdot {\bf x} \leq
  0, \, x_1 = 1, \, i = 1,2,\ldots, d \} $.  But by
  Lemma~\ref{halfspaces}, $Q[d]_{rev} = \emptyset$ which creates a
  contradiction.  Therefore, $x_j{\bf x}_{[d]}$ is not a minimal
  generator of $J$ and all preimages $P$ of $x_j$ have degree at most
  $d$.
\end{proof}

\noindent{\bf Example~\ref{example1} continued.} 
For $\sigma = \{1,4,13\}$, $\sigma_{out} = \{8,11,12\}$ which index 
the variables $h,k,l$. The minimal generators of $J$ that map to 
these variables under $\pi_{\sigma}$ are $hm, ak, al, ha, dml$. 

\vspace{.2cm}

Finally we consider the minimal generators of $J$ that do not project
under $\pi_{\sigma}$ to minimal generators of $J^{\sigma}$ for any
$\sigma \in \textup{max} \, \Delta_{\succ}$. Such generators may exist.

\vspace{.2cm}

\noindent{\bf Example~\ref{example1} continued.} 
Consider the minimal generator $gh$ of $J$. Then $gh = \\
\pi_{\{1,4,13\}}(gh) = \pi_{\{4,11,13\}}(gh) = \pi_{\{4,11,12\}}(gh) =
\pi_{\{11,12,13\}}(gh)$ is not a minimal generator of $J^{\sigma}$ for
any of the four facets $\sigma$ of $\Delta_{\succ}$.

\vspace{.2cm}

\begin{lemma} \label{lemma3ofdnormalbound}
  Let ${\bf x}^{\bf m}$ be a minimal generator of $J$ whose image
  under $\pi_{\sigma}$ is not a minimal generator of $J^{\sigma}$ for
  any facet $\sigma$ of $\Delta_{\succ}$. Then ${\bf x}^{\bf m}$ is a
  quadratic squarefree monomial.
\end{lemma}

\begin{proof}
  Let $\tau$ and $\tau'$ be facets of $\Delta_{\succ}$ and let $i \in
  \tau_{in}$, $j \in \tau_{in}'$ with $i, j \notin \tau_{in} \cap
  \tau_{in}'$. Then $x_i x_j$ is not covered by any standard pair of
  $J$ and hence lies in $J$. Since $\A$ is graded, $(x_i, \tau)$ and
  $(x_j, \tau')$ are standard pairs of $J$ which implies that $x_i
  x_j$ is a minimal generator of $J$.  Further, $\pi_{\sigma}(x_ix_j)$
  is not a minimal generator of $J^{\sigma}$ for any $\sigma \in
  \textup{max} \, \Delta_{\succ}$. We will prove that $L := \{ x_i x_j
  \, : \, i \in \tau_{in}, \, j \in \tau_{in}' \, and \, i, j \notin
  \tau_{in} \cap \tau_{in}' \}$ is precisely the set of minimal
  generators of $J$ that do not project under $\pi_{\sigma}$ to a
  minimal generator of $J^{\sigma}$ for a $\sigma \in \textup{max} \,
  \Delta_{\succ}$. This will prove the lemma.
  
  Suppose ${\bf x}^{\bf m}$ is a minimal generator of $J$ such that
  $\pi_{\sigma}({\bf x^m})$ is not a minimal generator of $J^{\sigma}$
  for any $\sigma \in \textup{max} \, \Delta_{\succ}$. Let $Y :=
  supp({\bf x}^{\bf m})$.
  
  {\em Case (i) $Y \subseteq \sigma_{in}$ for some $\sigma \in
    \textup{max} \, \Delta_{\succ}$}: Then ${\bf x}^{\bf m} \in
  N^{\sigma} = J^{\sigma} \cap {\mathbb K}[ x_j: j \in \sigma_{in}]$.
  Since ${\bf x^m}$ is not a minimal generator of $J^\sigma$ (and
  hence $N^{\sigma}$), some minimal generator of $N^{\sigma}$ properly
  divides ${\bf x^m}$. By Lemma~\ref{lemma1ofdnormalbound}, every
  minimal generator of $N^{\sigma}$ is a minimal generator of $J$
  which contradicts that ${\bf x^m}$ is a minimal generator of $J$.
 
  {\em Case (ii) $Y \subseteq \sigma$ for some $\sigma \in
    \textup{max} \, \Delta_{\succ}$}: Then ${\bf x^m}$ is covered by
  the standard pair $(1, \sigma)$ which contradicts that ${\bf x^m}
  \in J$.
    
  {\em Case (iii) $Y \subseteq \sigma \cup \sigma_{in}$ for some
    $\sigma \in \textup{max} \, \Delta_{\succ}$, with $Y \cap \sigma
    \neq \emptyset$ and $Y \cap \sigma_{in} \neq \emptyset$}: Write
  ${\bf x}^{\bf m} = {\bf x}^{{\bf m}^{\prime}} {\bf x}^{{\bf
      m}^{\prime \prime}} $ where $\emptyset \neq supp({\bf x}^{\bf
    m^{\prime}}) \subseteq \sigma$ and $\emptyset \neq supp({\bf
    x}^{{\bf m}^{\prime \prime}})\subseteq \sigma_{in}$. Then ${\bf
    x}^{{\bf m}^{\prime \prime}} \in N^\sigma$ all of whose minimal
  generators are minimal generators of $J$. This implies that a
  divisor of ${\bf x}^{{\bf m}^{\prime \prime}}$ is a minimal generator
  of $J$. Therefore, ${\bf x^m}$ is not a minimal generator of $J$, a
  contradiction.
    
    The above cases have shown that there is no single $\sigma \in
    \textup{max} \, \Delta_{\succ}$ such that $Y \subseteq \sigma \cup
    \sigma_{in}$. Therefore, there exists two distinct
    $\sigma, \tau \in \textup{max} \, \Delta_{\succ}$ and two indices
    $i, j \in Y$ such that $i \in \sigma \cup \sigma_{in} \cap
    \tau_{out}$ and $j \in \tau \cup \tau_{in} \cap \sigma_{out}$.
    
    {\em Case (a) $i \in \sigma$}: Since $j \in \sigma_{out}$, $x_i
    x_j$ is not covered by any standard pair of $J$ and so lies in
    $J$. Since $x_ix_j$ divides ${\bf x^m}$ and ${\bf x^m}$ is a
    minimal generator of $J$ it must be that ${\bf x^m} = x_ix_j$. But
    then $\pi_{\sigma}({\bf x^m}) = x_j$, $j \in \sigma_{out}$ is a
    minimal generator of $J^{\sigma}$ which contradicts our choice of
    ${\bf x^m}$.  Therefore this case cannot arise.
    
    {\em Case (b) $j \in \tau$}: By a symmetric argument to the
    previous, this cannot happen.
    
    Therefore, the only possibility is that $i \in \sigma_{in}$ and $j
    \in \tau_{in}$. Since $i \in \tau_{out}$ and $j \in \sigma_{out}$,
    $i, j \notin \sigma_{in} \cap \tau_{in}$. By the argument in the
    beginning of the proof, $x_ix_j$ is a minimal generator of $J$ and
    so ${\bf x^m} = x_ix_j$. Thus ${\bf x^m}$ lies in the set $L$ as
    claimed.
\end{proof}

\noindent {\em Proof of Theorem~\ref{dnormalbound}.}
Lemmas \ref{lemma1ofdnormalbound}, \,\ref{lemma2ofdnormalbound}
and \ref{lemma3ofdnormalbound} account for all minimal generators of
the initial ideal $J$ and show that they all have degree at most $d$.
Since $\A$ is graded, the reduced Gr\"obner basis of $I_{\A}$ with
initial ideal $J$ consists of homogeneous binomials. Hence these
binomials have degree at most $d$.~\hspace{\fill}$\Box$

\section{$\Delta$-normal and non-$\Delta$-normal families}
\label{families}

In this last section we construct non-trivial families of both
$\Delta$-normal and non-$\Delta$-normal configurations.  Recall that
any configuration $\mathcal A$ is always $\Delta$-normal with respect
to all its regular unimodular triangulations $\Delta$. Also, a
configuration $\mathcal A = \{{\bf a}_1, \ldots, {\bf a}_n\} \subset
\Z^d$ for which $cone(\mathcal A)$ is simplicial is $\Delta$-normal
with respect to its coarsest (regular) triangulation $\Delta = \{\{1,
\ldots, d\}\}$ if we assume that $cone(\mathcal A) = cone(\{{\bf a}_1,
\ldots, {\bf a}_d\})$. Call $\mathcal A$ {\em simplicial} if
$cone(\mathcal A)$ is simplicial. We construct families of
$\Delta$-normal configurations that are not simplicial and do not
admit regular unimodular triangulations. By computer search, Firla and
Ziegler \cite{f} found hundreds of normal simplicial configurations
$\mathcal A$ in $\N^4$ and $\N^5$ (in the course of writing \cite{fz})
that admit no unimodular triangulations. Our first result in this
section is a construction that extends a Firla-Ziegler configuration
to a family of {\em non-simplicial} $\Delta$-normal configurations ---
one in $\Z^d$ for each $d \geq 5$ --- without unimodular
triangulations. In the second part of this section we construct a
family of normal configurations (an $\mathcal A \subset \mathbb Z^d$
for each $d \geq 11$) that are not $\Delta$-normal for any regular
triangulation $\Delta$.

\subsection{$\Delta$-normal families from Firla-Ziegler
  configurations:}

Each Firla-Ziegler normal simplicial $\mathcal A \subset \mathbb N^4$
without unimodular triangulations is the Hilbert basis of the cone
generated by ${\bf e}_1, {\bf e}_2, {\bf e}_3$, the first three unit
vectors in $\mathbb R^4$, and a vector ${\bf v} \in \N^4$ of the form
${\bf v} := (a,b,c,d)^t$ with $0 < a < b < c < d$. In this subsection
we let $\mathcal A$ denote such a Firla-Ziegler configuration and let
$\mathcal A_{ext} = \{{\bf e}_1, {\bf e}_2, {\bf e}_3, {\bf v}\}$ be
the extreme rays of $cone(\mathcal A)$. By definition, $\mathcal A$ is
the unique minimal generating set of the semigroup $cone(\mathcal
A_{ext}) \cap \Z^4$ and $cone(\mathcal A) \subset \R^4_{\geq
  0}$.

\begin{lemma}
The vector ${\mathbf 1} := (1,1,1,1)^t$ is contained in $\mathcal
A$.
\end{lemma}

\begin{proof}
  Since ${\mathbf 1} = \frac{1}{d}{\bf v} + \frac{d-c}{d}{\bf e}_3 +
  \frac{d-b}{d}{\bf e}_2 + \frac{d-a}{d}{\bf e}_1$ and $a < b < c <
  d$, ${\mathbf 1} \in cone(\mathcal A_{ext}) \cap \Z^4 = \N \A$.  For
  every ${\bf p} = (p_1,p_2,p_3,p_4) \in cone(\mathcal A_{ext}) =
  cone(\A)$ with $p_4 > 0$, $p_i > 0$ for $i = 1, \ldots,4$ since the
  $\R_{\geq 0}$-linear combination of elements in $\mathcal A_{ext}$
  that expresses ${\bf p}$ as an element of $cone(\mathcal A_{ext})$
  must involve a positive multiple of ${\bf v}$. On the other hand,
  the $\N$-linear combination of elements in $\A$ that expresses
  ${\mathbf 1}$ as an element of $\mathbb N \mathcal A$ is the sum of
  distinct vectors in $\mathcal A \cap \{0,1\}^4$. At least one of
  these $0-1$ vectors --- say ${\bf w}$ --- has a positive last
  co-ordinate which implies that ${\bf w} = \mathbf 1$. Therefore,
  $\mathbf 1$ is in $\A$, the minimal Hilbert basis of $cone(\mathcal
  A_{ext}) \cap \mathbb Z^4$.
\end{proof}

\begin{example} \label{firstfz}
  The first Firla-Ziegler $\mathcal A$ in $\mathbb N^4$ has ${\bf v} =
  (1,2,3,5)$ and $\mathcal A$ consists of the columns of the matrix
  $${\bf A} = \left (
\begin{array}{cccccccc}
1 & 0 & 0 & 1 & 1 & 1 & 1 & 1 \\
0 & 1 & 0 & 2 & 1 & 1 & 2 & 2 \\
0 & 0 & 1 & 3 & 1 & 2 & 2 & 3 \\
0 & 0 & 0 & 5 & 1 & 2 & 3 & 4
\end{array}
\right ) .$$ The Hilbert basis of any rational polyhedral cone can
be computed using the software package Normaliz~\cite{bk}.
\end{example}

From a Firla-Ziegler $\mathcal A$ we will now recursively construct
configurations $\A^d$ for each $d \geq 5$ such that $\A^d \subset
\mathbb N^d$ is $\Delta$-normal, $cone(\A^d)$ is not simplicial and
$\A^d$ has no unimodular triangulations.  For each $d \geq 5$ let
${\bf p}_d = {\bf e}_1 + \cdots + {\bf e}_4 \in \Z^d$, ${\bf p}_d^+ =
{\bf p}_d + {\bf e}_d \in \Z^d$ and ${\bf p}_d^- = {\bf p}_d - {\bf
  e}_d \in \Z^d$.  Here ${\bf e}_d$ is the $d$-th unit vector in
$\R^d$.  Letting $\A^4 := \A$ (a Firla-Ziegler configuration in
$\mathbb N^4$) and $\A^4_{ext} := \A_{ext}$, recursively define
${\A^{d-1}}' := \{ ({\bf a},0) : {\bf a} \in \A^{d-1} \}$,
${\A_{ext}^{d-1}}' := \{ ({\bf a},0) : {\bf a} \in \A_{ext}^{d-1} \}$
and $\A^d := \{ {\bf p}_d^+ , {\bf p}_d^- \} \cup {\A^{d-1}}'$.  We
assume that ${\bf p}_d^+$ and ${\bf p}_d^-$ are always the first and
second elements of $\mathcal A^d$ and that $\sigma$ is the index set
of ${\A_{ext}^{d-1}}'$ in $\A^d$. Let $\sigma_1 = \{1\} \cup \sigma$
and $\sigma_2 = \{2\} \cup \sigma$.  Consider the triangulation
$\Delta^d$ of $cone(\mathcal A^d)$ consisting of the maximal subcones
$K_1 = cone(\A^d_{\sigma_1})$ and $K_2 = cone(\A^d_{\sigma_2})$.

\begin{lemma} \label{a5}
The configuration $\A^5$ has the following properties:
\begin{enumerate}
\item $\mathbb Z (\A^5 \cap K_1)  = \mathbb Z (\A^5 \cap K_2) = \mathbb Z^5$,
\item $\A^5$ is non-simplicial,
\item $\A^5$ is $\Delta^5$-normal, and
\item $\A^5$ admits no unimodular triangulations.
\end{enumerate}
\end{lemma}

\begin{proof}
\begin{enumerate}
\item Since ${\bf p}_5 = (1,1,1,1,0)$, ${\bf p}_5^+ = (1,1,1,1,1)$ and
  the first three unit vectors of $\mathbb R^5$ belong to $\A^5 \cap
  K_1$, it follows that all unit vectors of $\mathbb R^5$ lie in
  $\mathbb Z (\A^5 \cap K_1)$ which gives the result. Similarly,
  $\mathbb Z (\A^5 \cap K_2) = \mathbb Z^5$.
\item Since ${\bf p}_5$ lies in the interior of $cone(\A^5)$, the
  vectors ${\bf p}_5^+$ and ${\bf p}_5^-$ do not lie on a common facet
  of the cone. Hence $cone(\A^5)$ is a bipyramidal cone over
  $cone({\A^4}')$ with six extreme rays and is hence non-simplicial.
\item The triangulation $\Delta^5$ is the regular triangulation of
  $\A^5$ with respect to the weight vector ${\bf w} = {\bf e}_1 + {\bf
    e}_2$.  We first argue that $\A^5 \cap K_1$ is a minimal
  generating set of the semigroup $K_1 \cap \mathbb Z (\A^5 \cap K_1)
  \stackrel{(1)}{=} K_1 \cap \mathbb Z^5$. Suppose ${\bf q} = (q_1,
  \ldots, q_5) \in K_1 \cap \mathbb Z^5$. Since ${\bf p}_5^+$ is the
  unique generator of $K_1$ with a positive fifth co-ordinate, ${\bf
    q} = q_5{\bf p}_5^+ + {\bf q'}$ where ${\bf q}= (q_1-q_5, q_2-q_5,
  q_3-q_5, q_4-q_5,0)^t$ is the unique expression of ${\bf q}$ as an
  $\mathbb R_{\geq 0}$-combination of ${\bf p}_5^+$ and the other
  extreme rays of $K_1$.  Since $q'_5 = 0$, in fact, ${\bf q'} \in
  cone({\A^4}') \cap \mathbb Z^5 \stackrel{\ast}{=} \mathbb N {\A^4}'
  \subset \mathbb N (K_1 \cap \A^5)$ where the equality $(\ast)$
  follows from the normality of $\A^4$. This in turn implies that
  ${\bf q} = q_5{\bf p}_5^+ + {\bf q'} \in \mathbb N (K_1 \cap \A^5)$.
  (Note that $q_5 \in \N$.)  Similarly, $\A^5 \cap K_2$ is its own
  Hilbert basis.  Thus, $\A^5$ is $\Delta^5$-normal.
\item Suppose $T$ is a unimodular triangulation of $\A^5$ and $\tau$
  is a facet of $T$. Then by (1), $|det({\bf A}^5_{\tau})| \, = \, 1$ and
  $\{1,2\} \cap \tau \neq \emptyset$. If $\{1,2\} \subset \tau$, then
  $${\bf A}^5_{\tau} = \left ( \begin{array}{ccccc}
      1 & 1 & * & * & *\\
      1 & 1 & * & * & *\\
      1 & 1 & * & * & *\\
      1 & 1 & * & * & *\\
      1 & -1 & 0 & 0 & 0 \end{array} \right )$$
  which shows that
  $|det({\bf A}^5_{\tau})| \, \in \, 2 \Z $, a contradiction. Hence each
  maximal simplex in $T$ contains exactly one of ${\bf p}_5^+$ or
  ${\bf p}_5^-$ and $T$ induces a triangulation $T'$ of ${\A^4}'$.
  Since $T$ is unimodular, $T'$ gives a unimodular triangulation of
  $\A^4$ which is a contradiction as $\A^4$ has no unimodular
  triangulations.  Therefore, we conclude that $\A^5$ has no
  unimodular triangulations.
\end{enumerate}
\end{proof}

\begin{theorem} \label{babysfirstdeltanormalfamily}
For each $d \geq 5$, the configuration $\A^d$ has the following
properties:
\begin{enumerate}
\item $\mathbb Z (\A^d \cap K_1)  = \mathbb Z (\A^d \cap K_2) = \mathbb Z^d$,
\item $\A^d$ is non-simplicial,
\item $\A^d$ is $\Delta^d$-normal, and
\item $\A^d$ admits no unimodular triangulations.
\end{enumerate}
\end{theorem}

\begin{proof} This theorem is proved by induction using Lemma~\ref{a5}
  as the base step.
\begin{enumerate}
\item Suppose the result is true for $k = d-1$. Then it follows that
  $\mathbb Z {\A^{d-1}}'$ contains the first $d-1$ unit vectors of
  $\mathbb Z^d$ which are therefore also in $\mathbb Z (\A^d \cap
  K_1)$ and $\mathbb Z (\A^d \cap K_2)$. Since ${\bf p}_d^+ \in K_1
  \cap \A^d$ (and ${\bf p}_d^- \in K_2 \cap \A^d$), we also get that
  ${\bf e}_d \in \mathbb Z (\A^d \cap K_1)$ (and ${\bf e}_d \in
  \mathbb Z (\A^d \cap K_2)$). Hence $\mathbb Z (\A^d \cap K_1) =
  \mathbb Z (\A^d \cap K_2) = \mathbb Z^d$.
\item Assume by induction that $\A^{d-1}$ is non-simplicial and that
  ${\bf p}_{d-1}$ lies in the interior of $cone(\A^{d-1})$. Then,
  ${\bf p}_d$ lies in the interior of $cone(\A^d)$ and hence ${\bf
    p}_d^+$ and ${\bf p}_d^-$ do not lie on a common facet of
  $cone(\A^d)$.  This implies that $cone(\A^d) \subset \mathbb R^d$
  has exactly two more extreme rays than $cone(\A^{d-1}) \subset
  \mathbb R^{d-1}$ and hence is non-simplicial.
\item As in Lemma~\ref{a5}, $\Delta^d$ can be induced as the regular
  triangulation of $\A^d$ with respect to the weight vector ${\bf w} =
  {\bf e}_1 + {\bf e}_2$ for each $d \geq 5$. We assume by induction
  that $\A^{d-1}$ is $\Delta^{d-1}$-normal and hence normal. The
  arguments that $\A^d \cap K_1$ and $\A^d \cap K_2$ are minimal
  generating sets for $K_1 \cap \mathbb Z^d$ and $K_2 \cap \mathbb
  Z^d$ respectively follow from a straight generalization of the
  arguments in Lemma~\ref{a5}.
\item Again we assume by induction that $\A^{d-1}$ admits no
  unimodular triangulations. The rest of the argument is also a
  straight generalization of the arguments in Lemma~\ref{a5} (4).
\end{enumerate}

\end{proof}

We have thus produced non-simplicial $\Delta$-normal
configurations without unimodular triangulations in every
dimension beyond four, starting with a Firla-Ziegler $\mathcal A$
in $\mathbb N^4$.  The construction applies to all such Firla-Ziegler 
configurations.

\subsection{ Non $\Delta$-normal configurations from an example of
  Hibi and Ohsugi:}
\begin{figure}
\input{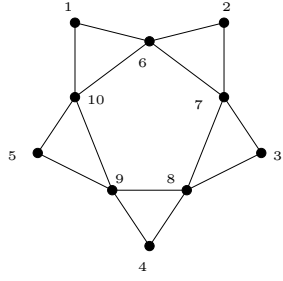}
\caption{The graph giving the Hibi-Ohsugi configuration} \label{ho}
\end{figure}

Consider the graph $G_{HO}$ shown in Figure~\ref{ho}.
In \cite{oh}  Hibi and Ohsugi showed that the graded normal
configuration
$$
\A_{HO} = \{ {\bf e}_1 + {\bf e}_i + {\bf e}_j \, : \, \{ i, j \}
\in E(G_{HO})\, , 1 \notin \{i, j \} \} \cup \{ {\bf e}_1 + {\bf e}_i
\, : \, \{ 1, i \} \in E(G_{HO}) \} $$
admits no regular unimodular
triangulations, although it does have non-regular unimodular
triangulations. Further, the $15$ points in $\A_{HO}$ are all extreme
points of the convex hull of $\A_{HO}$, denoted as $conv(\A_{HO})$.
The $(0,1)$-polytope $conv(\A_{HO}) \subset \{ {\bf x} \in \mathbb
R^{10} \, : \, x_1 = 1 \}$ is {\em empty} which means that it has no
lattice points other than its vertices. Thus $cone(\A_{HO})$ which is
a cone over $conv(\A_{HO})$ has $15$ extreme rays and is therefore
non-simplicial.

\begin{lemma} \label{emptynotdeltanormal}
  Let $\A \subset \Z^d$ be a normal graded non-simplicial
  configuration in $\{{\bf x} \in \mathbb R^d \, : \, x_1 = 1 \}$ such
  that $conv(\A)$ is empty. If $\A$ does not have a regular unimodular
  triangulation then $\A$ is not $\Delta$-normal for any regular
  triangulation $\Delta$.
\end{lemma}

\begin{proof} Without loss of generality, we can assume that $\Z \A =
  \Z^d$. By the hypothesis, every regular triangulation $\Delta$ of
  $\A$ has a maximal face $\sigma$ such that $|\,det({\bf A}_{\sigma})\,|
  \geq 2$. Thus the Hilbert basis of $cone(\A_{\sigma})$ contains at
  least one vector ${\bf q} \in \Z^d$ not in $\A_{\sigma}$. Since all
  vectors in $\A$ are extreme rays of $cone(\A)$, none of them lie in
  $cone(\A_{\sigma})$ unless they are in $\A_{\sigma}$. This implies
  that $\A_{\sigma} = cone(\A_{\sigma}) \cap \Z^d$ is not normal and
  hence $\A$ is not $\Delta$-normal.
\end{proof}

\begin{corollary}
  The Hibi-Ohsugi configuration $\A_{HO}$ is not $\Delta$-normal for
  any regular triangulation $\Delta$.
\end{corollary}

From $\A_{HO}$ we now recursively construct configurations $\A^d$
for each $d \geq 11$ such that $\A^d$ is normal and graded but not
$\Delta$-normal for any regular triangulation $\Delta$.  For each
$d \geq 11$ let ${\bf p}_d = {\bf e}_1 + {\bf e}_d \in \Z^d$.
Letting $\A^{10} := \A$, recursively define
$${\A^{d-1}}' := \left\{ \left( \begin{array}{c} {\bf a} \\ 0
    \end{array} \right) \,: \, {\bf a} \in \A^{d-1} \, \right\} \,\,
\textup{and} \,\, \A^d := \{{\bf p}_d\} \cup {\A^{d-1}}'.$$

\begin{theorem} \label{babysfirstnormalfamily}
  For each $d \geq 11$, the configuration $\A^d$ is normal and graded
  but not $\Delta$-normal for any regular triangulation $\Delta$.
\end{theorem}

\begin{proof}
  It suffices to show that $\A^d$ satisfies the conditions of
  Lemma~\ref{emptynotdeltanormal} for each $d$. For a given $d$,
  $\A^d$ is graded since it lies in $\{{\bf x} \in \mathbb R^d \, : \,
  x_1 = 1 \}$ and $conv(\A^d)$ is a $(0,1)$-polytope and hence empty.
  Further, $cone(\A^d)$ is non-simplicial as $cone(\A^{d-1})$ is
  non-simplicial for all $d \geq 11$.
  
  The configuration $\A^{11}$ is normal. To see this let ${\bf q} :=
  (q_1, \cdots, q_{11})^t \in cone(\A^{11}) \cap \Z \A^{11}$. Since
  the only extreme ray of $cone(\A^{11})$ with non-zero eleventh
  co-ordinate is ${\bf p}_{11}$, $q_{11} \leq \ q_1$. Further, ${\bf
    q} = q_{11}{\bf p}_{11} + {\bf q'}$ where ${\bf q'} = (q_1 -
  q_{11}, \, q_2, \cdots, q_{10}, \, 0)^t$ is the unique expression of
  ${\bf q}$ as an $\mathbb R_{\geq 0}$-combination of the extreme rays
  of $cone(\A^{11})$. The integral vector ${\bf q'}$ lies in $\mathbb
  N {\A^{10}}'$ since $\A^{10}$ is normal and hence it lies in
  $\mathbb N \A^{11}$. Thus ${\bf q} \in \mathbb N \A^{11}$.  By
  induction, it follows that $\A^d$ is normal for all $d \geq 11$.

  Suppose $\A^{11}$ had a regular unimodular triangulation. Then ${\bf 
    p}_{11}$ would be a vertex in every maximal face of this regular
  unimodular triangulation of $conv(\A)$. This in turn induces a
  regular unimodular triangulation in ${\A^{10}}'$ and hence in
  $\A^{10}$, a contradiction. Again, a straightforward inductive
  argument shows that $\A^d$ has no regular unimodular triangulation
  for all $d \geq 11$.
\end{proof}

\section*{Acknowledgments}
We thank Diane Maclagan for alerting us to the interpretation of
Simon's condition for a monomial ideal to be Cohen-Macaulay in her
paper \cite{ms} with Greg Smith. We also thank Winfried Bruns, Joseph
Gubeladze, Serkan Ho{\c s}ten, Francisco Santos and G{\"u}nter Ziegler
for discussions on constructing $\Delta$-normal configurations.

\end{document}